\documentclass{article}
\title{Incomplete ${\cal A}$-Hypergeometric Systems}
\author{Kenta Nishiyama\footnote{
  Graduate School of Information Science and Technology, Osaka
  University and JST CREST.
}
\; and Nobuki Takayama\footnote{
  Department of Mathematics, Kobe University and JST CREST.
  The second author is supported by KAKENHI 19204008.
}
}
\date{July 3, 2009, Revised September 1, 2011} 
\usepackage{amsmath,amscd,amssymb}

\newtheorem{theorem}{Theorem}

\newtheorem{remark}{Remark}
\newtheorem{example}{Example}
\newtheorem{definition}{Definition}
\newtheorem{corollary}{Corollary}

\def\pd#1{ \partial_{#1} }
\def\gausshg#1#2#3#4{
    F\left(
    \begin{matrix}
       #1 , & #2 \cr
              & #3
    \end{matrix} 
   \, ; #4
    \right)
}
\def\nquad#1{\count0=0%
          \loop \ifnum#1 > \count0 \quad \advance \count0 by 1 \repeat}

\begin{document}
\maketitle

\begin{abstract}
We generalize the theory of $\mathcal{A}$-hypergeometric systems
for incomplete functions.
We give a general algorithm on contiguity relations by utilizing Gr\"obner basis
as well as a detailed study on
incomplete Gauss and Lauricella hypergeometric functions.
\end{abstract}

\section{Introduction}
The incomplete gamma function, incomplete beta function and 
elliptic integrals are defined by definite integrals with parameters.
In modern theory of special functions, 
we regard hypergeometric functions as pairings of twisted cycles
and twisted cocycles \cite[Chap 2]{Aomoto-Kita}.
This is a fundamental fact and yields a lot of formulas of hypergeometric 
functions.
However, domains of integrations of incomplete gamma function, incomplete
beta function and elliptic integrals are not (twisted) cycles in general.
For example, the incomplete beta function
$ B(a,b; y) = \int_0^y t^{a-1} (1-t)^{b-1} dt $
is obtained by integration over $[0,y]$. 
When $y=1$ or $y=\infty$, it can be regarded as a (twisted) cycle, but if not, 
it cannot.
We will call this type of integral incomplete integral
in which the domain of integration contains parameters
and when parameters take special values the domain of integration 
can be regarded as a (twisted) cycle.

$\mathcal{A}$-hypergeometric systems (\cite{GZK}) 
contain a broad class of hypergeometric functions as their solutions and
it has been shown that $\mathcal{A}$ is a good framework to study hypergeometric functions
(see, e.g., \cite{SST}).
These solutions have integral representations
and domains of integrations are (twisted) cycles.
In recent studies of algebraic statistics, Gr\"obner basis and 
algebraic geometry play important roles. 
In these studies, we are now expected to study definite integrals
with parameters in statistics by algebraic methods. 
Some important integrals such as marginal likelihood functions 
in statistics, e.g., \cite{SSX},
look like $\mathcal{A}$-hypergeometric, but domains of integrations
are not (twisted) cycles.
We may regard  marginal likelihood functions as incomplete hypergeometric
functions.
Being motivated with these functions,
we want to generalize the algebraic 
theory of $\mathcal{A}$-hypergeometric systems
to that for incomplete integrals.

This paper is a first step toward this direction.
We give general discussions as well as detailed discussions
on $\Delta_1 \times \Delta_1$-incomplete hypergeometric functions,
which contain the incomplete beta function, the incomplete Gauss
and Lauricella functions, and elliptic integrals.
In the section \ref{section:general_definition}, 
we give a general definition
of incomplete $\mathcal{A}$-hypergeometric system
and prove that the solutions are holonomic functions.
In the section \ref{sec:Alg-contiguity},
we discuss an algorithm of deriving contiguity relations,
which is an incomplete version of the algorithm given in \cite{SST2}. 
The last section \ref{sec:1x1-ihg} consists of several subsections
which discuss about $\Delta_1 \times \Delta_1$-incomplete hypergeometric 
systems not only in the algebraic point of view but also in an analytic point 
of view.
We hope that our theorems and formulas are useful 
for numerical and asymptotic evaluations of 
incomplete $\mathcal{A}$-hypergeometric functions.

A definition of incomplete generalized hypergeometric functions is
given by Chardhry and Qadir \cite{chaudhry-qadir}.
A study of relations of their definition and ours is a future problem.

\section{General definition}  \label{section:general_definition}

Let $D$ be the Weyl algebra in $n$ variables.
A multi-valued holomorphic function $f$ defined on a Zariski open set
in ${\bf C}^n$ is called 
a {\it holonomic function}
if there exists a left ideal $I$ of $D$ such that
(1) $D/I$ is holonomic and (2) $I \bullet f = 0$.

We denote by $A=(a_{ij})$ a $d\times n$-matrix whose elements are integers.
We suppose that the set of the column vectors of $A$ spans ${\bf Z}^d$.

\begin{definition} \rm  \label{def:ihg}
We call the following system of differential equations $H_{A}(\beta,g)$
an {\it incomplete  $\mathcal{A}$-hypergeometric system}:
\begin{eqnarray*}
  \left( \sum_{j=1}^n a_{ij}  x_j \partial_{j} - \beta_i \right) \bullet f &=& g_i,
   \qquad(i = 1, \ldots, d)  \\
  \left( \prod_{i=1}^n \partial_{i}^{u_{i}} 
       - \prod_{j=1}^n \partial_{j}^{v_{j}} 
  \right) \bullet f &=& 0 
\end{eqnarray*}
$$\mbox{ with } u, v \in {\bf N}_0^{n} \mbox{ running over all $u, v$ such that }  Au =  Av.$$
Here, ${\bf N}_0 = \{ 0, 1, 2, \ldots \}$, and 
$\beta = (\beta_1, \ldots, \beta_d) \in {\bf C}^d$ are parameters
and $g=(g_1, \ldots, g_d)$ where $g_i$ are given holonomic functions
which may depend on parameters $\beta$.
We call solutions of the incomplete $\mathcal{A}$-hypergeometric system
{\it incomplete $\mathcal{A}$-hypergeometric functions}.
\end{definition}

Although we have introduced the incomplete $\mathcal{A}$-hypergeometric
system for arbitrary holonomic functions $g_i$,
$g_i$ are often also solutions of smaller incomplete
$\mathcal{A}$-hypergeometric system or well-known special functions
in interesting cases.

\begin{example} \rm
The {\it incomplete beta function} is defined as
$$ B(\alpha,\beta;y) = \int_0^y s^{\alpha-1} (1-s)^{\beta-1} ds $$
Replace $s=yt$. Then, we have
$B(\alpha,\beta;y) = y ^{\alpha} \int_0^1  t^{\alpha-1}
                        (1-yt)^{\beta-1} dt$.
Put $B(\alpha,\beta;y) = y^{\alpha} \tilde{B} (\beta-1,\alpha-1;1,-y)$
where 
$\tilde{B}(\gamma-1,\alpha-1;x_1,x_2)
=\int_0^1 t^{\alpha-1}(x_1+x_2t)^{\gamma-1} dt.$ 
The function $\tilde{B}$ is a solution of
an incomplete $\mathcal{A}$-hypergeometric system $H_A(\beta,g)$ for 
$A=\begin{pmatrix} 
    1 & 1 \cr
    0 & 1 \cr
   \end{pmatrix}
$,
$\beta=(\gamma-1,\alpha-1)$,
$g_1=0, g_2= (x_1+x_2)^{\beta-1}$.
Thus, the incomplete beta function can be expressed in terms of the
incomplete $\mathcal{A}$-hypergeometric function.
We will revisit this example in Example \ref{example:beta-contiguity}.
\end{example}

\begin{theorem}  \label{theorem:holonomic}
Solutions of the incomplete $\mathcal{A}$-hypergeometric system
are holonomic functions.
\end{theorem}

{\it Proof}.  
Since the $\mathcal{A}$-hypergeometric left ideal $H_A(\beta) \subset D$
is holonomic ideal,
the function $f$ satisfies the ordinary differential equation
of the form
$$\left(\sum_{k=0}^{r_i} a_k(\beta,x) \pd{i}^k \right) \bullet f
   = \sum_{i=1}^d \ell_i \bullet g_i
$$
where $\ell_i \in D$.
Since $g:=\sum_{i=1}^d \ell_i \bullet g_i$ is also a holonomic function,
there exists an ordinary differential operator such that
$ \left(\sum_{k=0}^{r'_i} b_k(\beta,x) \pd{i}^k \right) \bullet g = 0$.
Hence, $f$ is annihilated by an ordinary differential operator
with parameters of the order $r_i r'_i$ with respect to $x_i$
and then $f$ is annihilated by a zero dimensional ideal $J$
in ${\bf C}(x_1,\ldots, x_n) \langle \pd{i}, \ldots, \pd{n} \rangle$
\cite[p.33]{SST}.
Since $D/(J \cap D)$ is a left holonomic $D$-module 
(see, e.g.,  \cite[Th 2.4]{takayama-1992}, \cite[p.34]{SST}),
we have the conclusion.
We note that this construction can be done algorithmically
\cite{Oaku-Takayama-Walther}.
\bigbreak

The column vectors of the matrix $A$ are denoted by
$a_1, \ldots, a_d$.
We are interested in a probability distribution on ${\bf R}_{>0}^d = \{ t \}$ proportional to
$$ q(x,\beta;t)  = \exp\left(\sum_{i=1}^n x_i t^{a_i}\right) t^{-\beta-1} $$
where $\beta \in {\bf R}^d$ and $x=(x_1, \ldots, x_n)$ are parameters of the distribution.
The exponential family is an important class of probabilistic distributions.
The function $q$ defines a subclass of the exponential family.
The normalization integral (constant) for the subclass 
satisfies our incomplete $\mathcal{A}$-hypergeometric system.

\begin{theorem} \label{theorem:intrep}
Let $P$ be a polytope in ${\bf R}_{> 0}^d$. 
Put 
$(dt)_k = (-1)^{k-1}dt_1 \wedge \cdots \wedge dt_{k-1} \wedge dt_{k+1} \wedge 
  \cdots \wedge dt_d $
and
$$ g_k = \int_{\partial P} q(x,\beta;t) t_k (dt)_k. $$
Then, the functions $g_k$ are holonomic and the function
$$ f(x) = \int_P q(x,\beta;t) dt $$
satisfies the incomplete $\mathcal{A}$-hypergeometric equations.
\end{theorem}

{\it Proof}\/.
Let us show that $g_k$ is holonomic.
Since $P$ is a polytope, we may assume that the boundary $\partial P$ is
a finite union of simplices.
Since the sum of finite holonomic functions is also holonomic,
we may prove it when $\partial P$ is parametrized as
$t_i = \sum_{j=1}^{d-1} b_{ij} s_j$, $b_{ij} \in {\bf R}$,
$s \in \Delta = \{ s \mid 0 \leq s_j \leq 1, (j=1, \ldots, d-1),
                          1-s_1-\cdots -s_{d-1} \geq 0 \}$.
Then, the differential form
$(dt)_k$ can be written as $b ds_1 \cdots ds_{d-1}$ where $b$ is a constant.
The integral  $g_k$ is equal to 
$ \int_{\Delta} q(x,\beta; t(s)) t_k(s) b ds $.
Let $h(s)$ be the product of the Heaviside functions
$H(s_1), \ldots, H(s_{d-1})$ and $ H(1-\sum_{i=1}^{d-1} s_i)$.
Then, we have
$ \int_{\Delta} q(x,\beta; t(s)) t_k(s) b ds 
= \int_{{\bf R}^{d-1}} q(x,\beta; t(s)) t_k(s) b h(s) ds$.
Since the integrand is a holonomic distribution, the integral on ${\bf R}^{d-1}$
is a holonomic function.
Thus, we have proved that $g_k$ is a holonomic function.

We will prove that $f$ is a solution of the incomplete $\mathcal{A}$-hypergeometric equations.
Since we have $\left(\pd{}^u - \pd{}^v \right) \bullet q(x,\beta;t) = 0$, 
we have $\left(\pd{}^u - \pd{}^v \right) \bullet f = 0$.
Act $\sum_{j=1}^n a_{kj} x_j \pd{j}$ to $f$.
Then, we obtain
$ \int_P  \frac{\partial \exp(\sum_{i=1}^n x_i t^{a_i})}
               { \partial t_k} t_k t^{-\beta-1} dt$.
Since 
$d ( q(x,\beta;t) t_k (dt)_k ) 
= \frac{\partial q(x,\beta;t) t_k}
       {\partial t_k} 
  dt
$,
we have, by the Stokes theorem,
\begin{eqnarray*}
 \int_{\partial P} q(x,\beta;t) t_k (dt)_k
&=& \int_P d (q(x,\beta;t) t_k (dt)_k) \\
&=& \int_P  \frac{\partial \exp(\sum_{i=1}^n x_i t^{a_i})}
               { \partial t_k} t_k t^{-\beta-1} dt
- \beta_k \int_P  q(x,\beta;t) dt.
\end{eqnarray*}
This shows that the function $f$ satisfies the $k$-equation of the incomplete hypergeometric equations.
\bigbreak

It is known that 
integrals of products of polynomials over cycles
satisfy ${\cal A}$-hypergeometric systems.
We have an analogous result for incomplete hypergeometric systems.
Let us state the result.
Let $ A_1 = (a_1, \ldots, a_{n_1}), \ldots,
A_\ell = (a_{n_{\ell-1}+1}, \ldots, a_{n_\ell})$ \,
($a_i\in {\bf Z}^d$) be $\ell$ matrices and
$$
A = \left(
\begin{array}{cccccccccccc}
1 & \cdots & 1 & 0 & \cdots & 0 & & & & 0 & \cdots & 0 \\
0 & \cdots & 0 & 1 & \cdots & 1 & & & & 0 & \cdots & 0 \\
0 & \cdots & 0 & 0 & \cdots & 0 & & & & 0 & \cdots & 0 \\
\cdot & \cdots & \cdot & \cdot & \cdots & \cdot & & & & \cdot & \cdots & \cdot \\
\cdot & \cdots & \cdot & \cdot & \cdots & \cdot & & \cdots & & \cdot & \cdots & \cdot \\
\cdot & \cdots & \cdot & \cdot & \cdots & \cdot & & & & \cdot & \cdots & \cdot \\
0 & \cdots & 0 & 0 & \cdots & 0 & & & & 1 & \cdots & 1 \\
a_1 & \cdots & a_{n_1}& a_{n_1+1}& \cdots & a_{n_2}& & & & a_{n_{\ell-1}+1}& \cdots & a_{n_\ell}
\end{array}
\right)
$$
be a $(\ell+d) \times n_\ell$ matrix.
To each matrix $A_j$ we associate a generic polynomial with that support:
$$ f_j(x,t) = \sum_{i=n_{j-1}+1}^{n_j} x_i t^{a_i}. $$
For complex numbers
$\alpha_1, \ldots, \alpha_\ell$, and $\gamma_1, \ldots, \gamma_d$,
we consider the polynomial
$$q(x, \alpha, \gamma, t) 
= t^{\gamma} \prod_{j=1}^\ell f_j(x,t)^{\alpha_j}.$$

\begin{theorem}
Let $P$ be a polytope in ${\bf R}^d \setminus \{ t \prod f_j = 0 \}$. 
Put 
$$ g_k =
\left\{
\begin{array}{cl}
 0, &\quad {\text if \; } 1 \leq k \leq \ell; \\
 \int_{\partial P} q(x, \alpha, \gamma; t) t_{k-\ell} (dt)_{k-\ell}, &\quad
  {\text if \; } \ell+1 \leq k \leq \ell+d.
\end{array}
\right.$$
Then, the functions $g_k$ are holonomic and the function
$$ f(x) = \int_P q(x, \alpha, \gamma; t) dt $$
satisfies the incomplete $\mathcal{A}$-hypergeometric system
$H_A(\beta, g)$ where
$\beta=(\alpha_1, \ldots, \alpha_\ell, -\gamma_1-1, \ldots, -\gamma_d-1)$.
\end{theorem}

The proof is analogous to that of Theorem \ref{theorem:intrep}.


\section{Algorithms deriving contiguity relations} \label{sec:Alg-contiguity}
In this section, we will give algorithms to obtain contiguity relations
for incomplete $\mathcal{A}$-hypergeometric functions under the condition that
for all $i$ $(1 \leq i \leq d)$ there exists a constant $c_k \in {\bf C} \setminus \{0\}$
such that 
\begin{equation} \label{eq:condition_contiguity}
\pd{k}\bullet g_i(\beta) = c_k^{-1} g_i(\beta-a_k).
\end{equation}
This condition holds in many interesting cases.
We mean by a contiguity relation, a relation among two incomplete 
$\mathcal{A}$-hypergeometric functions $\Phi(\beta; x)$ and $\Phi(\beta^\prime; x)$
where $\beta - \beta^{\prime} \in {\bf Z}^d$
and $\Phi(\beta;x)$ is a solution of $H_A(\beta,g)$.
Under the condition (\ref{eq:condition_contiguity}), we note that
$$c_k^{-1} \pd{k}\bullet \Phi(\beta; x) =:  \Phi(\beta-a_k; x) $$
is a solution of $H_A(\beta-a_k,g)$.
In other words, the operator $\pd{k}$ gives a contiguity relation 
for the parameter shift $-a_k$.
This contiguity relation follows from
\begin{align*}
 & \pd{k} \left (\sum_{j=1}^n a_{ij}  x_j \partial_{j} - \beta_i \right ) \bullet f \\
 =& \left (\sum_{j=1}^n a_{ij}  x_j \partial_{j} - (\beta_i - a_{ik}) \right ) \pd{k} \bullet f
 = \pd{k}\bullet g_i(\beta)
\end{align*}
and the condition (\ref{eq:condition_contiguity}).

If we can find the inverse operator for $\pd{k}$, then it means that
we find ``generators'' of the all contiguity relations because the matrix $A$
is the full rank.
In \cite{SST2} and \cite{SST}, they give algorithms to find the inverse operator
for $\pd{k}$ in case of $\mathcal{A}$-hypergeometric system.
We utilize their algorithms to find the inverse operator.
Let us recall Definition \ref{def:ihg}. 
Let $h_i$ $(1 \leq i \leq d)$ be $\sum_{j=1}^n a_{ij}  x_j \partial_{j} - \beta_i$ 
respectively 
and $h_i$ $(d+1 \leq i \leq s)$ be generators of the toric ideal $I_A$ of the form
$\prod_{i=1}^n \partial_{i}^{u_{i}} - \prod_{j=1}^n \partial_{j}^{v_{j}}$.
It follows from \cite{SST2} and \cite{SST} that there exist $r, r_i \in D(\beta)$
such that
$$ r \pd{k} + \sum_{i=1}^s r_i h_i = 1.$$
The operators $r, r_i$ can be obtained by the syzygy computation by Gr\"obner basis.
Applying the operator of the left hand side to $\Phi(\beta; x)$,
we obtain the contiguity relation 
$$ c_k r \bullet \Phi(\beta-a_k; x) + \sum_{i=1}^{d} r_i \bullet g_i = \Phi(\beta; x)$$
which contains the functions $g_i$.

Our second algorithm gives contiguity relations which do not contain
the functions $g_i$.
Suppose that we are given operators ${\tilde h}_i$ such that ${\tilde h}_i \bullet g_i = 0$ and
operators
$ \pd{k} $, ${\tilde h}_i h_i$ ($ 1 \leq i \leq d$),
and $h_i$ ($d+1 \leq i \leq s$) generate a trivial ideal in $D(\beta)$.
Then, 
we can construct operators $r, r_i \in D(\beta)$
such that
$$ r \pd{k} + \sum_{i=1}^d r_i {\tilde h}_i h_i  + \sum_{i=d+1}^s r_i h_i = 1$$ by the Gr\"obner basis method.
The operator $c_k r$ is the inverse of $\pd{k}$.

We will apply these algorithms to obtain
a complete list of contiguity relations of 
$\Delta_1 \times \Delta_1$-hypergeometric functions
in Section \ref{sec:1x1-ihg}.

\section{Incomplete $\Delta_1 \times \Delta_1$-hypergeometric functions} \label{sec:1x1-ihg}

In the previous sections, we have given a general discussion on incomplete hypergeometric
systems.
An important example is the incomplete beta function,
which is defined as an integral of a product of two power functions.
It will be natural to consider a product of $3$ power functions
and regard it as an incomplete Gauss hypergeometric function.
We will give a detailed study on this function from our point of view.

We assume $0 < a < b $ for simplicity.
We consider the integral
\begin{equation}  \label{eq:intighg}
 \int_a^b t^\gamma (x_{11}+x_{21}t)^{\alpha_1}(x_{12}+x_{22}t)^{\alpha_2} dt
\end{equation}
for $x_{ij} > 0$ and ${\rm Re}\, \gamma , {\rm Re}\, \alpha_i > -1$.
The integral and its analytic continuations satisfy the following incomplete 
$\mathcal{A}$-hypergeometric system.
\begin{equation} \label{eq:inhomogeneous}
\left\{
\begin{array}{ll}
  \left (\theta_{11} \theta_{22} - \frac{x_{11} x_{22}}{x_{21} x_{12}} \theta_{21} \theta_{12} \right) \bullet f & = 0 \\
  (\theta_{11} + \theta_{21} - \alpha_1) \bullet f & =  0 \\
  (\theta_{12} + \theta_{22} - \alpha_2) \bullet f & =  0 \\
  (\theta_{21} + \theta_{22} + \gamma + 1) \bullet f & = [g(t,x)]_{t=a}^{t=b}
\end{array}
\right.
\end{equation}
Here, $g(t,x) = t^{\gamma+1}(x_{11}+x_{21}t)^{\alpha_1}(x_{12}+x_{22}t)^{\alpha_2}$
and $\theta_{ij} = x_{ij} \pd{ij}$.
This fact can be shown by exchanging the integral and 
differentiations (see, e.g., \cite[p.221]{SST}).
When $[g(t,x)]_{t=a}^{t=b}=0$, our system is essentially 
the Gauss hypergeometric equation (see, e.g., \cite[Chap 1]{SST}).

Since the matrix 
$A=\begin{pmatrix} 
    1 & 1 & 0 & 0 \cr
    0 & 0 & 1 & 1 \cr
    0 & 1 & 0 & 1 \cr
   \end{pmatrix}
$
associated to the system can be regarded as
$\Delta_1 \times \Delta_1$ ($1$-simplex times $1$-simplex), 
we call it
the incomplete $\Delta_1 \times \Delta_1$-hypergeometric system and
its solutions are 
incomplete $\Delta_1 \times \Delta_1$-hypergeometric functions.

\begin{example}  \rm \label{ex:i-elliptic}
The {\it incomplete elliptic integral of the first kind} is defined as
$$F(z;k) = \int_0^z \frac{dx}{\sqrt{(1-x^2)(1-k^2x^2)}}.$$
Replacing $x^2$ by $z^2t$, we obtain
$$F(z;k) = \frac{1}{2} z \int_0^1 t^{-\frac{1}{2}}(1-z^2t)^{-\frac{1}{2}}(1-k^2z^2t)^{-\frac{1}{2}} dt,$$
which agrees with $\frac{z}{2} \times \mbox{(\ref{eq:intighg})}$
with $\gamma = \alpha_1 = \alpha_2 = - \frac{1}{2}$,
$a=0, b=1$,
$x_{11} = 1, x_{21}= -z^2, x_{12}=1, x_{22}=-k^2 z^2$,
$ [g]_{t=a}^{t=b} = g_{|_{t=1}}$.
Thus, the incomplete elliptic integral of the first kind
can be regarded as a solution restricted to a subvariety of incomplete 
$\mathcal{A}$-hypergeometric system $H_A(\beta,g)$ for 
$A=\begin{pmatrix} 
    1 & 1 & 0 & 0 \cr
    0 & 0 & 1 & 1 \cr
    0 & 1 & 0 & 1 \cr
   \end{pmatrix}
$,
$\beta = (-\frac{1}{2},-\frac{1}{2},-\frac{1}{2})$.
We will revisit this example in Example \ref{example:expression_F_1}.
\end{example}

In order to make a rigorous discussion, we need to specify branches 
of multi-valued functions appearing in our discussion. 
In the sequel, $z^\alpha$ denotes the unique analytic continuation of
the function $z^\alpha$ defined on $z>0$ to the upper and
the lower half plane as long as we make no annotation.
\begin{remark} \rm \label{remark:branch}
Under this definition,
we have
$(zw)^\alpha = e^{2 \pi \alpha} z^\alpha w ^\alpha$
for ${\rm Im}\, z <0, {\rm Im}\, w < 0$ and ${\rm Im}\, zw > 0$,
$(zw)^\alpha = e^{- 2 \pi \alpha} z^\alpha w ^\alpha$
for ${\rm Im}\, z >0 , {\rm Im}\, w > 0$ and ${\rm Im}\, zw < 0$,
and
$(zw)^\alpha = z^\alpha w ^\alpha$
for other cases.
\end{remark}

We take real numbers
$ x^*_{ij} > 0 $ such that
$$0 < \frac{x^*_{12}}{x^*_{22}}  < \frac{x^*_{11}}{x^*_{21}} 
  < a < b 
  < \frac{x^*_{21}}{x^*_{11}} < \frac{x^*_{22}}{x^*_{12}} 
$$

We consider the simply connected domain
define by 
\begin{eqnarray*} 
&& 
 (-1)^{d_1} {\rm Im}\, x_{11} > 0, 
 (-1)^{d_2} {\rm Im}\, x_{21} > 0, \  
 (-1)^{d_3} {\rm Im}\, x_{12} > 0, \  
 (-1)^{d_4} {\rm Im}\, x_{22} > 0, \   \\
&&
 (-1)^{d_5} {\rm Im}\, x_{21}/x_{11} > 0, \  
 (-1)^{d_6} {\rm Im}\, x_{22}/x_{12} > 0, \  
 (-1)^{d_7} {\rm Im}\, \frac{x_{21} x_{12}}{x_{11} x_{22}} > 0
\end{eqnarray*}
Here, $d_i$ takes the values $0$ or $1$.
Since we assume  $a, b \in {\bf R}$
and the singular locus of the homogeneous $\mathcal{A}$-hypergeometric 
system is 
$x_{11} x_{21} x_{12} x_{22} (x_{11} x_{22} - x_{21} x_{12}) = 0$,
the solutions of our system are holomorphic on each of these domains.
We denote by $D_d$ where $d=(d_1, \ldots, d_7) \in \{0, 1\}^7$
the simply connected domain standing for $d$.

We define the four domains as follows
\begin{eqnarray*}
D_{12}^{11} &=& \{ x_{ij} \,|\, 
   |x_{21}b /x_{11}| < 1, 
   |x_{21}a /x_{11}| < 1, 
   |x_{22}b /x_{12}| < 1, 
   |x_{22}a /x_{12}| < 1 
 \} \\
D_{22}^{11} &=& \{ x_{ij} \,|\, 
   |x_{21}b /x_{11}| < 1, 
   |x_{21}a /x_{11}| < 1, 
   |x_{12} /(x_{22}b)| < 1, 
   |x_{12} /(x_{22}a)| < 1 
 \} \\
D_{12}^{21} &=& \{ x_{ij} \,|\, 
   |x_{11} /(x_{21} b)| < 1, 
   |x_{11} /(x_{21} a)| < 1, 
   |x_{22}b /x_{12}| < 1, 
   |x_{22}a /x_{12}| < 1 
 \} \\
D_{22}^{21} &=& \{ x_{ij} \,|\, 
   |x_{11} /(x_{21} b)| < 1, 
   |x_{11} /(x_{21} a)| < 1, 
   |x_{12} /(x_{22}b)| < 1, 
   |x_{12} /(x_{22}a)| < 1 
 \} \\
\end{eqnarray*}

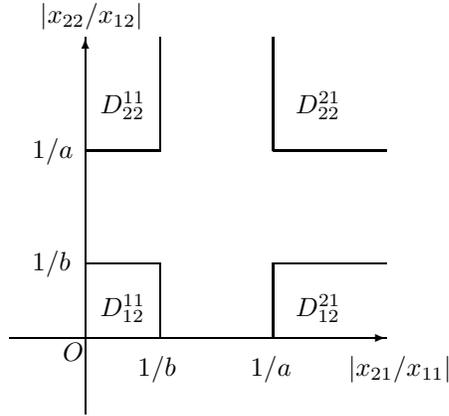
\begin{figure}[htbp]
 \begin{center}
\setlength\unitlength{1truecm}
\begin{picture}(5,5)(0,0)
  \put(0.4,5.2){{$\vert x_{22}/x_{12} \vert$}}
  \put(4.5,0.5){{$\vert x_{21}/x_{11}\vert$}}
  \put(0.7,0,7){{$O$}}
  \put(3.2,0.5){{$1/a$}}
  \put(1.7,0.5){{$1/b$}}
  \put(0.3,1,9){{$1/b$}}
  \put(0.3,3.4){{$1/a$}}
  \put(1.2,1.3){{$D_{12}^{11}$}}
  \put(1.2,4.0){{$D_{22}^{11}$}}
  \put(3.8,1.3){{$D_{12}^{21}$}}
  \put(3.8,4.0){{$D_{22}^{21}$}}
  \put(1,0){\vector(0,1){5}}
  \put(0,1){\vector(1,0){5}}
  \put(2,1){\line(0,1){1}}
  \put(1,2){\line(1,0){1}}
  \put(2,3.5){\line(0,1){1.5}}
  \put(1,3.5){\line(1,0){1}}
  \put(3.5,2){\line(1,0){1.5}}
  \put(3.5,1){\line(0,1){1}}
  \put(3.5,3.5){\line(1,0){1.5}}
  \put(3.5,3.5){\line(0,1){1.5}}
\end{picture}
\end{center}
\caption{Domains $D_{ij}^{k\ell}$ on $(|x_{21}/x_{11}|,|x_{22}/x_{12}|)$-plane}
\label{fig:domain}
\end{figure}

The point $(x^*_{ij})$ belongs to the last domain $D_{22}^{21}$.
It is easy to see that each of these domains
and $D_d$ has an open intersection since $a,b \in {\bf R}$.

\begin{remark} \rm \label{remark:branchOfg}
We note that
$$(x_{1i}+x_{2i}t)^{\alpha_i}
 =x_{1i}^{\alpha_i} \left(1 + \frac{x_{2i} t}{x_{1i}} \right)^{\alpha_i}
 =(x_{2i}t)^{\alpha_i} \left( \frac{x_{1i}}{x_{2i} t} + 1 \right)^{\alpha_i}
$$
for $x_{ij} > 0$ and $t>0$.
This relation will be used to specify branches of $[g(t,x)]_{t=a}^{t=b}$.
For example,
when $(x_{ij}) \in D_{22}^{21}$, we regard
$[g(t,x)]_{t=a}^{t=b}$ as 
$$\left[t^{\gamma+1}(x_{21}t)^{\alpha_1}\left(\frac{x_{11}}{x_{21} t} + 1\right)^{\alpha_1}
                   (x_{22}t)^{\alpha_2}\left(\frac{x_{12}}{x_{22} t} + 1\right)^{\alpha_2}
 \right]_{t=a}^{t=b}
$$
of which series expansion converges on $D_{22}^{21}$.
Since the domain $D_{22}^{21}$ has an open intersection with $D_d$,
the function $[g(t,x)]_{t=a}^{t=b}$ has a unique analytic continuation to $D_d$. 
\end{remark}

\subsection{Homogeneous system}

As we have proved in Theorem \ref{theorem:holonomic},
solutions of incomplete $\mathcal{A}$-hypergeometric systems
are holonomic functions.
The advantage of this point of view is that
we can apply some algorithms for holonomic systems to study solutions
of our system.
For example, we can apply the algorithm given in the Chapter 2 of \cite{SST}
to find candidates of series solutions.
Holonomic systems which annihilate these functions can be obtained 
in an algorithmic way.
However, outputs by the algorithm are sometimes tedious.

In the case of $\Delta_1 \times \Delta_1$-hypergeometric system,
solutions satisfy the following relatively simple holonomic system.
\begin{equation} \label{eq:homogeneous}
\left\{
\begin{array}{ll}
  \left (\theta_{11} \theta_{22} - \frac{x_{11} x_{22}}{x_{21} x_{12}} \theta_{21} \theta_{12} \right) \bullet f & = 0 \\
  (\theta_{11} + \theta_{21} - \alpha_1) \bullet f & =  0 \\
  (\theta_{12} + \theta_{22} - \alpha_2) \bullet f & =  0 \\
  (\pd{22}-a\pd{12})(\pd{21}-b\pd{11})(\theta_{21} + \theta_{22} + \gamma + 1) \bullet f & = 0
\end{array}
\right.
\end{equation}

\subsection{Contiguity relations}

We will derive contiguity relations of incomplete $\Delta_1 \times \Delta_1$-hypergeometric functions by applying our two algorithms.
In this section, we put $\beta = -\gamma-1$ to make formulas of contiguity
relations of the incomplete $\Delta_1 \times \Delta_1$-hypergeometric function
simpler forms. We put 
$$\Phi(\alpha_1,\alpha_2,\beta;x) = \int_a^b t^{-\beta-1}(x_{11}+x_{21}t)^{\alpha_1}(x_{12}+x_{22}t)^{\alpha_2} dt.$$

\begin{theorem} \label{theorem:contiguity-a_i}
The incomplete $\Delta_1 \times \Delta_1$-hypergeometric function 
$\Phi(\alpha_1,\alpha_2,\beta;x)$ satisfies the following contiguity relations.
\begin{itemize}
\item Shifts with respect to $a_1=(1,0,0)$
\begin{align*}
 S(\alpha_1,\alpha_2,\beta;-a_1) \Phi(\alpha_1,\alpha_2,\beta)
 &= \alpha_1 \Phi(\alpha_1-1,\alpha_2,\beta)\\
 \hat{S}(\alpha_1-1,\alpha_2,\beta;+a_1) \Phi(\alpha_1-1,\alpha_2,\beta)
 &= (\alpha_1+\alpha_2-\beta) \Phi(\alpha_1,\alpha_2,\beta) - [g(t,x)]_{t=a}^{t=b}\\
 S(\alpha_1-1,\alpha_2,\beta;+a_1) \Phi(\alpha_1-1,\alpha_2,\beta)
 &= \alpha_2(\alpha_1+\alpha_2-\beta) \Phi(\alpha_1,\alpha_2,\beta)
\end{align*}
where

\noindent
$S(\alpha_1,\alpha_2,\beta;-a_1) = \pd{11}$,\\
$\hat{S}(\alpha_1-1,\alpha_2,\beta;+a_1) 
  = (x_{21}x_{12}-x_{11}x_{22})\pd{22}+(\alpha_1+\alpha_2)x_{11}$,\\
$S(\alpha_1-1,\alpha_2,\beta;+a_1) 
  = x_{21}x_{22} (((a+b)(x_{21}\partial_{21}+x_{22}\partial_{22}+1-\beta)
    +x_{11}\partial_{21}+x_{12}\partial_{22})\partial_{22}
    +ab(\beta \partial_{12}-x_{22}\partial_{12}\partial_{22}-x_{21}\partial_{11}\partial_{22}))
    +(\alpha_1+\alpha_2-\beta)(x_{21}x_{12}\partial_{22}+\alpha_2 x_{11})$

\item Shifts with respect to $a_2=(1,0,1)$
\begin{align*}
 S(\alpha_1,\alpha_2,\beta;-a_2) \Phi(\alpha_1,\alpha_2,\beta)
 &= \alpha_1 \Phi(\alpha_1-1,\alpha_2,\beta-1)\\
 \hat{S}(\alpha_1-1,\alpha_2,\beta-1;+a_2) \Phi(\alpha_1-1,\alpha_2,\beta-1)
 &= \beta \Phi(\alpha_1,\alpha_2,\beta) + [g(t,x)]_{t=a}^{t=b}\\
 S(\alpha_1-1,\alpha_2,\beta-1;+a_2) \Phi(\alpha_1-1,\alpha_2,\beta-1)
 &= ab \alpha_2 \beta \Phi(\alpha_1,\alpha_2,\beta)
\end{align*}
where

\noindent
$S(\alpha_1,\alpha_2,\beta;-a_2) = \pd{21}$,\\
$\hat{S}(\alpha_1-1,\alpha_2,\beta-1;+a_2) 
  = x_{11}x_{22}\pd{12}+x_{21}x_{22}\pd{22}+\alpha_1 x_{21}$,\\
$S(\alpha_1-1,\alpha_2,\beta-1;+a_2)
  = (a+b) x_{11} (\alpha_2 x_{21} \partial_{21}-x_{21} x_{22} \partial_{21} \partial_{22}-x_{22}^2 \partial_{22}^2+(\beta+\alpha_2-2) x_{22} \partial_{22}-\alpha_2 (\beta-1))
     +a b (x_{11} x_{22} (x_{22} \partial_{22}-(\alpha_2+\beta-1)) \partial_{12}+(\alpha_2 \partial_{21}-x_{22} \partial_{21} \partial_{22}) x_{21}^2+x_{21} ((\alpha_1-1) x_{22} \partial_{22}-\alpha_2 (\alpha_1+\beta-1)))
     -x_{11} (x_{12} x_{22} \partial_{22}^2+(\alpha_1-\beta+1) x_{12} \partial_{22}+x_{11} x_{22} \partial_{21} \partial_{22}-\alpha_2 x_{11} \partial_{21})$


\item Contiguity relations with respect to $a_3=(0,1,0)$ are
obtained from those with respect to $a_1$ by the permutations 
$\alpha_1 \leftrightarrow \alpha_2$, $x_{i1} \leftrightarrow x_{i2}$,
$\pd{i1} \leftrightarrow \pd{i2}$.

\item Contiguity relations with respect to $a_4=(0,1,1)$ are
obtained from those with respect to $a_2$ by the same permutations 
as above.
\end{itemize}
\end{theorem}

{\it Proof}.
Since the function $[g(t,x)]_{t=a}^{t=b}$ satisfies the condition
(\ref{eq:condition_contiguity}), we can apply the first algorithm for 
$A = (a_1,a_2,a_3,a_4) = 
   \begin{pmatrix} 
    1 & 1 & 0 & 0 \cr
    0 & 0 & 1 & 1 \cr
    0 & 1 & 0 & 1 \cr
   \end{pmatrix}
$
to obtain the contiguity relations containing the function $[g(t,x)]_{t=a}^{t=b}$.
Contiguity relations which do not contain the function $[g(t,x)]_{t=a}^{t=b}$
is obtained from (\ref{eq:homogeneous}).
The generation condition of the trivial ideal generated by 1 is checked 
for (\ref{eq:homogeneous}) by a computer and then we can apply the second
algorithm in section \ref{sec:Alg-contiguity}.
\bigbreak

Theorem \ref{theorem:contiguity-a_i} gives contiguity relations for 
$e_1=(1,0,0)$, $e_2=(0,1,0)$, but it does not give those for $e_3=(0,0,1)$.
The set of vectors $\{e_1, e_2, e_3 \}$ is the standard basis of ${\bf Z}^3$.
The contiguity relations for $e_3$ can be obtained from Theorem 
\ref{theorem:contiguity-a_i} as follows.

\begin{corollary} \label{cor:contiguity-e_3}
The incomplete $\Delta_1 \times \Delta_1$-hypergeometric function 
$\Phi(\alpha_1,\alpha_2,\beta;x)$ satisfies the following contiguity relations.
\begin{itemize}
\item Shifts with respect to $e_3=(0,0,1)$ 
\begin{align*}
 S(\alpha_1,\alpha_2,\beta+1;-e_3) \Phi(\alpha_1,\alpha_2,\beta+1)
 &= \alpha_1 \alpha_2 (\alpha_1+\alpha_2-\beta) \Phi(\alpha_1,\alpha_2,\beta) \\
 S(\alpha_1,\alpha_2,\beta-1;+e_3) \Phi(\alpha_1,\alpha_2,\beta-1)
 &= ab\alpha_1 \alpha_2 \beta \Phi(\alpha_1,\alpha_2,\beta) \\
 \hat{S}(\alpha_1,\alpha_2,\beta+1;-e_3) \Phi(\alpha_1,\alpha_2,\beta+1)
 &= \alpha_1(\alpha_1+\alpha_2-\beta) \Phi(\alpha_1,\alpha_2,\beta)
 -\alpha_1 [g(t,x)]_{t=a}^{t=b} \\
 \hat{S}(\alpha_1,\alpha_2,\beta-1;+e_3) \Phi(\alpha_1,\alpha_2,\beta-1)
 &= \alpha_1 \beta \Phi(\alpha_1,\alpha_2,\beta) 
 +\alpha_1 [g(t,x)]_{t=a}^{t=b}
\end{align*}
where

\noindent
$S(\alpha_1,\alpha_2,\beta+1;-e_3) = S(\alpha_1-1,\alpha_2,\beta;+a_1) S(\alpha_1,\alpha_2,\beta+1;-a_2)$,\\
$S(\alpha_1,\alpha_2,\beta-1;+e_3) = S(\alpha_1-1,\alpha_2,\beta-1;+a_2) S(\alpha_1,\alpha_2,\beta-1;-a_1)$,\\
$\hat{S}(\alpha_1,\alpha_2,\beta+1;-e_3) = \hat{S}(\alpha_1-1,\alpha_2,\beta;+a_1) S(\alpha_1,\alpha_2,\beta+1;-a_2)$,\\
$\hat{S}(\alpha_1,\alpha_2,\beta-1;+e_3) = \hat{S}(\alpha_1-1,\alpha_2,\beta-1;+a_2) S(\alpha_1,\alpha_2,\beta-1;-a_1)$.\\
\end{itemize}
\end{corollary}

\begin{example} \rm  \label{example:beta-contiguity}
Let $a=0, b=1$. We consider the following degenerated incomplete $\Delta_1 \times \Delta_1$-hypergeometric function.
$$\Psi(\alpha_1,\beta;x) = \int_0^1 t^{-\beta-1} (x_{11}+x_{21}t)^{\alpha_1} dt$$
Then the last contiguity relation of Corollary \ref{cor:contiguity-e_3} for this function is
$$ x_{21} \pd{11} \Psi(\alpha_1,\beta-1;x) = \beta \Psi(\alpha_1,\beta;x) + (x_{11}+x_{21})^{\alpha_1}. $$
Multiplying both sides by $x_{11}$
and by using the relation of the incomplete $\mathcal{A}$-hypergeometric system:
$$x_{11}\pd{11}
 \Psi(\alpha_1,\beta-1;x)=(\alpha_1-\beta+1)\Psi(\alpha_1,\beta-1;x)-(x_{11}+x_{21})^{\alpha_1}, $$
 we have
$$ (\alpha_1-\beta+1) x_{21} \Psi(\alpha_1,\beta-1;x) 
  = \beta x_{11} \Psi(\alpha_1,\beta;x) + (x_{11}+x_{21})^{\alpha_1+1}. $$
Put $x_{11}=1, x_{21}=-y$ and
replace $\beta$ by $-\alpha$, $\alpha_1$ by $\beta-1$, we have
$$ (\alpha+\beta) (-y) \Psi(\beta-1,-\alpha-1;y)
  = -\alpha \Psi(\beta-1,-\alpha;y) + (1-y)^{\beta}.$$
Multiplying both sides by $-y^\alpha$, we obtain
$$ (\alpha+\beta) B(\alpha+1,\beta;y)
= \alpha B(\alpha,\beta;y) - y^{\alpha} (1-y)^{\beta}.$$
This is a well-known relation of the incomplete beta function.
\end{example}

\subsection{Series solutions}

We define the following 4 series.
\begin{align*}
  f_{12}^{11} & = x_{11}^{\alpha_1}x_{12}^{\alpha_2} \sum_{k,m \geq 0} 
         \frac{(-1)^{k+m}}{\gamma+k+m+1} 
         \cdot \frac{(-\alpha_1)_k(-\alpha_2)_m}{(1)_k(1)_m} \nonumber \\
         &\nquad{12} \cdot (b^{\gamma+k+m+1} - a^{\gamma+k+m+1})
         \left ( \frac{x_{21}}{x_{11}} \right )^{k}
         \left ( \frac{x_{22}}{x_{12}} \right )^{m} \\
  f_{22}^{11} & = x_{11}^{\alpha_1}x_{22}^{\alpha_2} \sum_{k,m \geq 0} 
         \frac{(-1)^{k+m}}{\gamma+\alpha_2+k-m+1} 
         \cdot \frac{(-\alpha_1)_k(-\alpha_2)_m}{(1)_k(1)_m} \nonumber \\
         &\nquad{12} \cdot (b^{\gamma+\alpha_2+k-m+1} - a^{\gamma+\alpha_2+k-m+1})
         \left ( \frac{x_{21}}{x_{11}} \right )^{k}
         \left ( \frac{x_{12}}{x_{22}} \right )^{m} \\
  f_{12}^{21} & = x_{21}^{\alpha_1}x_{12}^{\alpha_2} \sum_{k,m \geq 0} 
         \frac{(-1)^{k+m}}{\gamma+\alpha_1-k+m+1} 
         \cdot \frac{(-\alpha_1)_k(-\alpha_2)_m}{(1)_k(1)_m} \nonumber \\
         &\nquad{12} \cdot (b^{\gamma+\alpha_1-k+m+1} - a^{\gamma+\alpha_1-k+m+1})
         \left ( \frac{x_{11}}{x_{21}} \right )^{k}
         \left ( \frac{x_{22}}{x_{12}} \right )^{m} \\
  f_{22}^{21} & = x_{21}^{\alpha_1}x_{22}^{\alpha_2} \sum_{k,m \geq 0} 
         \frac{(-1)^{k+m}}{\gamma+\alpha_1+\alpha_2-k-m+1} 
         \cdot \frac{(-\alpha_1)_k(-\alpha_2)_m}{(1)_k(1)_m} \nonumber \\
         &\nquad{12} \cdot (b^{\gamma+\alpha_1+\alpha_2-k-m+1} - a^{\gamma+\alpha_1+\alpha_2-k-m+1})
         \left ( \frac{x_{11}}{x_{21}} \right )^{k}
         \left ( \frac{x_{12}}{x_{22}} \right )^{m}
\end{align*}
Here, $(a)_n$ is the Pochhammer symbol $ (a)_n := a(a+1) \cdots (a+n-1)$.
\begin{theorem} \label{theorem:series_sol}
We assume that $\gamma \not \in {\bf Z}$,
$\gamma+\alpha_i \not \in {\bf Z}$,
$\gamma+\alpha_1 + \alpha_2 \not \in {\bf Z}$.
\begin{enumerate}
\item The series $f_{ij}^{k\ell}$ converges on the domain $D_{ij}^{k\ell}$
and has a unique analytic continuation to $D_d$.
Here, $D_{ij}^{k\ell}$ is defined in the beginning of this section.
\item The function $f_{ij}^{k\ell}$ defined on $D_d$ as above satisfies
the incomplete $\Delta_1 \times \Delta_1$-hypergeometric system
for the branch of $[g(t,x)]_{t=a}^{t=b}$ given in the Remark \ref{remark:branchOfg}.
\item
$f_{12}^{11}$ can be expressed in terms of the Appell function $F_1$ as
\begin{eqnarray*}
 f_{12}^{11}
&=& x_{11}^{\alpha_1} x_{12}^{\alpha_2}
\left(
 \frac{b^{\gamma+1}}{\gamma+1}
 F_1\left(\gamma+1,-\alpha_1, -\alpha_2,\gamma+2,
          \frac{-x_{21}b}{x_{11}},
          \frac{-x_{22}b}{x_{12}}
    \right) \right. \\
& & \quad\quad - \left.
 \frac{a^{\gamma+1}}{\gamma+1}
 F_1\left(\gamma+1,-\alpha_1, -\alpha_2,\gamma+2,
          \frac{-x_{21}a}{x_{11}},
          \frac{-x_{22}a}{x_{12}}
    \right) \right).
\end{eqnarray*}
\end{enumerate}
\end{theorem}

{\it Proof}. 
The item 1 is proved by utilizing majorant series.  
There exists a constant $C$ such that
\begin{eqnarray*}
 &&
 C
 \left(\sum_{k=0}^\infty \frac{|(-\alpha_1)_k|}{k!}
       \left| \frac{b x_{21}}{x_{11}}   \right|^k 
 \right)   
 \left(\sum_{m=0}^\infty \frac{|(-\alpha_2)_m|}{m!}
       \left| \frac{b x_{22}}{x_{12}}   \right|^m 
 \right) \\
 & + &
  C
 \left(\sum_{k=0}^\infty \frac{|(-\alpha_1)_k|}{k!}
       \left| \frac{a x_{21}}{x_{11}}   \right|^k 
 \right)   
 \left(\sum_{m=0}^\infty \frac{|(-\alpha_2)_m|}{m!}
       \left| \frac{a x_{22}}{x_{12}}   \right|^m 
 \right) 
\end{eqnarray*}
is a majorant series of $f_{12}^{11}$.
Other cases can be shown analogously.

The item 2 is proved by
applying the algorithm to find series solutions for (\ref{eq:homogeneous})
given in the Chapter 2 of \cite{SST}.
The Gr\"obner cone consists of $8$ maximal dimensional cones.
After constructing series solutions of the homogeneous system
(\ref{eq:homogeneous}),
we check if they satisfy the inhomogeneous system
(\ref{eq:inhomogeneous}) and we find these four solutions.

The item 3 can be proved by utilizing the relation
$ \frac{1}{\gamma+k+m+1} 
  = \frac{1}{\gamma+1} \frac{(\gamma+1)_{k+m}}{(\gamma+2)_{k+m}}
$.
\bigbreak

\begin{example} \rm \label{example:expression_F_1}
As we have seen in Example \ref{ex:i-elliptic},
the incomplete elliptic integral of the first kind can be regarded
as incomplete $\Delta_1 \times \Delta_1$-hypergeometric function.
Let us apply Theorem \ref{theorem:series_sol} to obtain an expression
of the incomplete elliptic integral
in terms of the Appell function $F_1$.

Put $x_{11}=1,x_{21}=-z^2,x_{12}=1,x_{22}=-k^2z^2$ and 
$\alpha_1=\alpha_2=\gamma=-\frac{1}{2}$, $a=0,b=1$.
Then we have
\begin{align*}
F(z;k) &= \frac{1}{2}z \cdot \frac{1}{-\frac{1}{2}+1}
 F_1\left(-\frac{1}{2}+1,\frac{1}{2}, \frac{1}{2},-\frac{1}{2}+2;
          \frac{z^2}{1},\frac{k^2z^2}{1}
    \right)\\
    &= z  F_1\left(\frac{1}{2},\frac{1}{2}, \frac{1}{2},\frac{3}{2};
          z^2,k^2z^2
    \right).
\end{align*}
This expression of the incomplete elliptic integral 
seems to be well-known \cite{WFS}.
\end{example}

\begin{remark} \rm
The common refinement of the Gr\"obner fan of $H_A(\beta)$ and
that of ${\rm Ann}\, [g(t,x)]_{t=a}^{t=b}$
is a set of natural domains of definitions of series solutions in this case. 
\end{remark}

\subsection{Connection formulas}

Connection formulas for the Gauss hypergeometric functions
are given on the upper half plane and on the lower half plane.
We will give connection formulas of our series solutions
in an analogous way. 

The domain of convergence of our series solution $f_{ij}^{k\ell}$
and $D_d$ 
is non-empty and open set for any $d$, 
then there exists a unique analytic continuation of the series
$f_{ij}^{k\ell}$ to the domain $D_d$.
We will give connection formulas among our $4$ series solutions on $D_d$.

\begin{theorem}
We suppose $0 < a < b $ and exponents $\alpha_1, \alpha_2, \gamma$
are generic.
\begin{enumerate}
\item
\begin{eqnarray*}
  f_{12}^{11} &=&  e^{2\pi i \alpha_1} f_{12}^{21}
  \quad \mbox{ on $D_{(1,1,*,*,0,*,*)}$} \\
  f_{12}^{11} &=&  e^{-2\pi i \alpha_1} f_{12}^{21}
  \quad \mbox{ on $D_{(0,0,*,*,1,*,*)}$} \\
  f_{12}^{11} &=&  f_{12}^{21}
  \quad \mbox{ on other $D_d$'s}
\end{eqnarray*}
\item  
\begin{eqnarray*}
  f_{12}^{11} &=&  e^{2\pi i \alpha_2} f_{22}^{11}
  \quad \mbox{ on $D_{(*,*,0,1,*,0,*)}$} \\
  f_{12}^{11} &=&  e^{-2\pi i \alpha_2} f_{22}^{11}
  \quad \mbox{ on $D_{(*,*,1,0,*,1,*)}$} \\
  f_{12}^{11} &=&  f_{22}^{11}
  \quad \mbox{ on other $D_d$'s}
\end{eqnarray*}
\item
\begin{eqnarray*}
  f_{12}^{21} &=&  e^{2\pi i \alpha_2} f_{22}^{21}
  \quad \mbox{ on $D_{(*,*,0,1,*,0,*)}$} \\
  f_{12}^{21} &=&  e^{-2\pi i \alpha_2} f_{22}^{21}
  \quad \mbox{ on $D_{(*,*,1,0,*,1,*)}$} \\
  f_{12}^{21} &=&  f_{22}^{21}
  \quad \mbox{ on other $D_d$'s}
\end{eqnarray*}
\end{enumerate}
\end{theorem}

Intuitively speaking, the series $f_{ij}^{k\ell}$ are different expansions of
the same integral (\ref{eq:intighg})
in different domains and hence they will agree with some adjustments
of constant factor as in the Theorem.
Here, we will give a proof without using the integral representation.
The advantage of this discussion is that we can avoid topological discussions 
about choices of branches of the integrand.
Analogous discussion is used to study global behavior of solutions
of the Euler-Darboux equation \cite{takayama-euler-darboux}. 

{\it Proof}.
We note
$ \frac{1}{\gamma+m+k+1}
= \frac{1}{\gamma+m+1}
   \frac{(\gamma+m+1)_k}{(\gamma+m+2)_k}
$.
Then, the series $f_{12}^{11}$ can be expressed as a superposition of
contiguous family of Gauss hypergeometric functions as follows.
\begin{eqnarray}
&&x_{11}^{\alpha_1} x_{12}^{\alpha_2} \left(
\sum_{m=0}^\infty \left(\frac{-x_{22}}{x_{12}}\right)^m
                                \frac{b^{\gamma+m+1} (-\alpha_2)_m}{(\gamma+m+1) m!}
                                \gausshg{-\alpha_1}{\gamma+m+1}{\gamma+m+2}{\frac{-x_{21} b}{x_{11}}} 
                                \right.  \nonumber \\
&&\quad
- \left.
\sum_{m=0}^\infty \left(\frac{-x_{22}}{x_{12}}\right)^m
                                \frac{a^{\gamma+m+1} (-\alpha_2)_m}{(\gamma+m+1) m!}
                                \gausshg{-\alpha_1}{\gamma+m+1}{\gamma+m+2}{\frac{-x_{21} a}{x_{11}}} 
\right)   \label{eq:superposition}
\end{eqnarray}
The Gauss hypergeometric function
has the unique analytic continuation to
${\rm Im}\, x_{21}/x_{11} > 0$ and
${\rm Im}\, x_{21}/x_{11} < 0$.
We replace the Gauss hypergeometric functions in (\ref{eq:superposition})
with their series expansions around $x_{21}/x_{11} = \infty$.
In other words, we make replacements by using the connection formula
of the Gauss hypergeometric function in (\ref{eq:superposition}).
For the first hypergeometric function in (\ref{eq:superposition}),
we utilize
\begin{eqnarray*}
&& \gausshg{-\alpha_1}{\gamma+m+1}{\gamma+m+2}{\frac{-x_{21} b}{x_{11}}} \\
&=& \frac{\Gamma(\gamma+m+2) \Gamma(\gamma+m+1+\alpha_1)}
                {\Gamma(\gamma+m+1)\Gamma(\gamma+m+2+\alpha_1)}
         \left( \frac{x_{21} b}{x_{11}} \right)^{\alpha_1}
         \gausshg{-\alpha_1}{1-\alpha_1-\gamma-m-2}{1-\alpha_1-\gamma-m-1}{\frac{-x_{11}}{x_{21} b}}\\
 & & + 
  \frac{\Gamma(\gamma+m+2) \Gamma(-\alpha_1-\gamma-m-1)}
                {\Gamma(-\alpha_1)\Gamma(\gamma+m+2-\gamma-m-1)}
         \left( \frac{x_{21} b}{x_{11}} \right)^{-\gamma-m-1} \\
 & & \quad\quad\quad\quad \cdot
         \gausshg{\gamma+m+1}{1+\gamma+m+1-\gamma-m-2}{1+\gamma+m+1+\alpha_1}{\frac{-x_{11}}{x_{21} b}}\\
&=& \frac{\gamma+m+1}
                {\gamma+\alpha_1+m+1}
         \left( \frac{x_{21} b}{x_{11}} \right)^{\alpha_1}
         \gausshg{-\alpha_1}{1-\alpha_1-\gamma-m-2}{1-\alpha_1-\gamma-m-1}{\frac{-x_{11}}{x_{21} b}}\\
 & & + 
  \frac{\Gamma(\gamma+m+2) \Gamma(-\alpha_1-\gamma-m-1)}
                {\Gamma(-\alpha_1)}
         \left( \frac{x_{21} b}{x_{11}} \right)^{-\gamma-m-1} \\
\end{eqnarray*}
and the analogous formula for the second Gauss hypergeometric function in (\ref{eq:superposition}).
The terms obtained from the second terms of the connection formulas of the Gauss hypergeometric
functions are canceled and we obtain
\begin{eqnarray*}
&&x_{11}^{\alpha_1} x_{12}^{\alpha_2} \left(\frac{x_{21}}{x_{11}} \right)^{\alpha_1} 
\left(
\sum_{m=0}^\infty \left(\frac{-x_{22}}{x_{12}}\right)^m
                                \frac{b^{\alpha_1+\gamma+m+1} (-\alpha_2)_m}{(\gamma+\alpha_1+m+1) m!}
                                \gausshg{-\alpha_1}{1-\alpha_1-\gamma-m-2}{-\alpha_1-\gamma-m}{\frac{-x_{11} }{x_{12}b}} 
                                \right.   \\
&&\quad
- \left.
\sum_{m=0}^\infty \left(\frac{-x_{22}}{x_{12}}\right)^m
                                \frac{a^{\alpha_2+\gamma+m+1} (-\alpha_2)_m}{(\gamma+\alpha_1+m+1) m!}
                                \gausshg{-\alpha_1}{1-\alpha_1-\gamma-m-2}{-\alpha_1-\gamma-m}{\frac{-x_{11} }{x_{12}a}} 
\right).
\end{eqnarray*}
Expanding the Gauss hypergeometric functions,
we see that the above sum equals to
$ x_{11}^{\alpha_1} x_{12}^{\alpha_2} (x_{21}/x_{11})^{\alpha_1}
  x_{21}^{-\alpha_1} x_{12}^{-\alpha_2} f_{12}^{21} $.
Applying the formulas in Remark \ref{remark:branch}, 
we obtain the first result 1.
Other cases can be obtained analogously.

\subsection{Monodromy formula}
It is well known that the monodromy representation of 
the (complete) $\Delta_1 \times \Delta_n$-hypergeometric system 
can be understood as a $1$-cocycle of the Braid group $B_{n+1}$
(see, e.g., \cite{AGH}).
Any solution of the incomplete $\Delta_1 \times \Delta_n$-hypergeometric system 
is written as a sum of a solution of the complete system and 
a constant multiple of $f_{ij}^{k\ell}$.
Then, in order to study analytic continuation of the incomplete system,
we may study analytic continuation (monodromy) of the function $f_{ij}^{k\ell}$.
We only give formula for $f_{12}^{11}$.
Formulas for other $f_{ij}^{k\ell}$ can be obtained analogously by symmetry.
In order to give formulas, we define
\begin{align*}
f_{12}^{11}(p,q;x) 
  &= x_{11}^{\alpha_1}x_{12}^{\alpha_2} \sum_{k,m \geq 0} 
         \frac{(-1)^{k+m}}{\gamma+k+m+1} 
         \cdot \frac{(-\alpha_1)_k(-\alpha_2)_m}{(1)_k(1)_m} \nonumber \\
         &\nquad{12} \cdot (qb^{\gamma+k+m+1} - pa^{\gamma+k+m+1})
         \left ( \frac{x_{21}}{x_{11}} \right )^{k}
         \left ( \frac{x_{22}}{x_{12}} \right )^{m}, \\
\tilde{f}(x) &= x_{11}^{\alpha_1}x_{12}^{\alpha_2}
    \left ( -\frac{x_{11}}{x_{21}} \right )^{\gamma+1}
    \gausshg{-\alpha_2}{\gamma+1}{\gamma+\alpha_1+2}
    {\frac{x_{11}x_{22}}{x_{12}x_{21}} }.
\end{align*}
We note that $\tilde{f}(x)$ is a solution of the homogeneous system $H_A(\beta)$.

\begin{theorem}
We fix $x_{12},x_{21},x_{22}$ to real numbers for simplicity and regarded the function 
as a function in one variable $x_{11}$.
Let $\gamma_a$ be a path which encircles the point $-ax_{21}$ 
in the positive direction and
$\gamma_b$ be a path which encircles the point $-bx_{21}$ 
in the positive direction.
We also suppose that exponents are generic.
The analytic continuations of $f_{12}^{11}$ along $\gamma_a$ and $\gamma_b$ are
\begin{eqnarray*}
f_{12}^{11}(1,1;x) &\looparrowright_{\gamma_b^*} & f_{12}^{11}(1,e^{2 \pi i \alpha_1};x) 
                        + \frac{-2\pi i e^{\pi i (\alpha_1+1)}}{\Gamma(-\alpha_1)}\tilde{f}(x),\\
f_{12}^{11}(1,1;x) &\looparrowright_{\gamma_a^*} & f_{12}^{11}(e^{2 \pi i \alpha_1},1;x) 
                        - \frac{-2\pi i e^{\pi i (\alpha_1+1)}}{\Gamma(-\alpha_1)}\tilde{f}(x).
\end{eqnarray*}
\end{theorem}
{\it Proof}.
We replace the Gauss hypergeometric function in (\ref{eq:superposition})
with analytic continuations of them;
we utilize the following formula of the analytic continuation of the Gauss
hypergeometric function $F(a,b,c;x)$ along a path which encircles $x=1$ positively.
$$F(a,b,c;x) \looparrowright (1-A)F(a,b,c;x) + B x^{1-c}F(a-c+1,b-c+1,2-c;x)$$
Here, we put
$$
A=\frac{(1-e^{-2 \pi ia})(1-e^{-2 \pi ib})}{1-e^{-2 \pi ic}},
B=\frac{2 \pi i}{1-c} \cdot \frac{\Gamma(c)^2}{\Gamma(a)\Gamma(b)\Gamma(c-a)\Gamma(c-b)} \cdot e^{\pi i (c-a-b)}.
$$
We note that 
\begin{align*}
1-A &= 1- \frac{(1-e^{-2 \pi i(-\alpha_1)})(1-e^{-2 \pi i(\gamma+m+1)})}{1-e^{-2 \pi i(\gamma+m+2)}} \\
  &= 1-(1-e^{2 \pi i \alpha_1})\\
  &= e^{2 \pi i \alpha_1}.
\end{align*}
Then, we obtain the first term $f_{12}^{11}(1,e^{2 \pi i \alpha_1};x)$
of the right hand side of the first formula by the replacement of the type
$(1-A)F(a,b,c;x)$.
We make replacements of the type $x^{1-c}F(a-c+1,b-c+1,2-c;x)$.
Then, we have 
$$\gausshg{-\alpha_1-\gamma-m-1}{0}{-\gamma-m}{-\frac{x_{21}b}{x_{11}}} = 1.$$
Therefore, we have
\begin{eqnarray*}
& & x_{11}^{\alpha_1}x_{12}^{\alpha_2} \sum_{m=0}^{\infty}
 \left ( -\frac{x_{12}}{x_{22}} \right )^m
 \frac{b^{\gamma+m+1} (-\alpha_2)_m}{(\gamma+m+1)m!} \\
& & \nquad{2}\frac{2\pi i}{-(\gamma+m+1)} \cdot
 \frac{\Gamma(\gamma+m+2)^2}{\Gamma(-\alpha_1)\Gamma(\gamma+m+1)\Gamma(\gamma+\alpha_1+m+2)\Gamma(1)}
 \cdot e^{\pi i (\alpha_1+1)} \left ( -\frac{x_{21}b}{x_{11}} \right )^{-(\gamma+m+1)} \\
&& = x_{11}^{\alpha_1}x_{12}^{\alpha_2} \sum_{m=0}^{\infty}
 \left ( -\frac{x_{12}}{x_{22}} \right )^m
 \frac{b^{\gamma+m+1} (-\alpha_2)_m}{(\gamma+m+1)m!} \\
& & \nquad{2}\frac{2\pi i}{-(\gamma+m+1)} \cdot
 \frac{(\gamma+m+1)^2\Gamma(\gamma+m+1)^2}{\Gamma(-\alpha_1)\Gamma(\gamma+m+1)\Gamma(\gamma+\alpha_1+m+2)}
 \cdot e^{\pi i (\alpha_1+1)} b^{-(\gamma+m+1)} \left ( -\frac{x_{11}}{x_{21}} \right )^{\gamma+m+1} \\
&& = \frac{-2\pi i e^{\pi i (\alpha_1+1)}}{\Gamma(-\alpha_1)}
  \left ( -\frac{x_{11}}{x_{21}} \right )^{\gamma+1}
 x_{11}^{\alpha_1}x_{12}^{\alpha_2} \sum_{m=0}^{\infty}
 \left ( -\frac{x_{12}}{x_{22}} \right )^m
 \frac{(-\alpha_2)_m}{m!} \cdot
 \frac{\Gamma(\gamma+m+1)}{\Gamma(\gamma+\alpha_1+m+2)}
 \left ( -\frac{x_{11}}{x_{21}} \right )^{m}\\
&& = \frac{-2\pi i e^{\pi i (\alpha_1+1)}}{\Gamma(-\alpha_1)}
  \tilde{f}(x)
\end{eqnarray*}
which is the second term of the right hand side of the first formula.
The second formula in the Theorem is obtained analogously by exchanging
the role of $b$ and $a$.

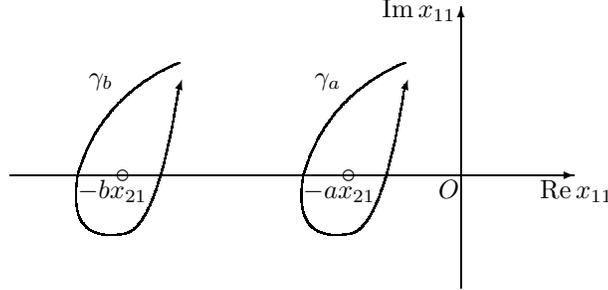
\begin{figure}[htbp]
 \begin{center}
\setlength\unitlength{1.5truecm}
\begin{picture}(5,2.3)(0,0)
  \put(3,1){\circle{0.1}}
  \put(1,1){\circle{0.1}}
  \put(3.8,0.8){$O$}
  \put(4.7,0.8){${\rm Re}\,x_{11}$}
  \put(3.3,2.4){${\rm Im}\,x_{11}$}
  \put(2.6,0.8){$-ax_{21}$}
  \put(0.6,0.8){$-bx_{21}$}
  \put(2.7,1.8){$\gamma_a$}
  \put(0.7,1.8){$\gamma_b$}
  \put(3.52,1.82){\vector(1,3){0.01}}
  \put(1.52,1.82){\vector(1,3){0.01}}
  \put(4,0){\vector(0,1){2.5}}
  \put(0,1){\vector(1,0){5}}
  \qbezier(3.5,2.0)(2.8,1.7)(2.6,1)
  \qbezier(2.6,1)(2.5,0.4)(3,0.48)
  \qbezier(3,0.48)(3.3,0.5)(3.5,1.8)
  \qbezier(1.5,2.0)(0.8,1.7)(0.6,1)
  \qbezier(0.6,1)(0.5,0.4)(1,0.48)
  \qbezier(1,0.48)(1.3,0.5)(1.5,1.8)
\end{picture}
\end{center}
\caption{Paths $\gamma_b$ and $\gamma_a$}
\label{fig:path}
\end{figure}

\begin{remark} \rm
The function $[g(t,x)]_{t=a}^{t=b}$ is analytically continued as follows.
\begin{eqnarray*}
{[g(t,x)]}_{t=a}^{t=b} &\looparrowright_{\gamma_b^*} & e^{2 \pi i \alpha_1} g(b,x)- g(a,x)\\
{[g(t,x)]}_{t=a}^{t=b} &\looparrowright_{\gamma_a^*} & g(b,x)- e^{2 \pi i \alpha_1}g(a,x) 
\end{eqnarray*}
\end{remark}


\begin{thebibliography}{99}
\bibitem{Aomoto-Kita}
K.~Aomoto, M.~Kita, 
{\it Theory of Hypergeometric Functions}.
Springer, Tokyo, 1994 (in Japanese).
%
\bibitem{chaudhry-qadir}
M.A.~Chaudhry, Asghar Qadir,
Incomplete Exponential and Hypergeometric Functions
with Applications to the Non Central $\chi^2$-Distribution.
Communications in Statistics -- Theory and Models {\bf 34}
(2005),  525--535.
%
\bibitem{GZK} 
I.M.~Gel'fand, A.V.~Zelevinsky, M.M.~Kapranov, 
Hypergeometric functions and toral manifolds.
Functional Analysis and its Applications {\bf 23} (1989), 94--106.
%
\bibitem{AGH}
R.P.~Holzapfel, A.M.~Uludag, M.~Yoshida (Editors),
{\it Arithmetic and Geometry around Hypergeometric Functions},
Progress in Mathematics 260, Birkh\"auser, 2007.
%
\bibitem{SSX}
Shaowei Lin, B.~Sturmfels, Zhiqiang Xu, 
Marginal Likelihood Integrals for Mixtures of Independence Models.
Journal of Machine Learning Research {\bf 10} (2009), 1611--1631.
%
\bibitem{Oaku-Takayama-Walther}
T.~Oaku, N.~Takayama, U.~Walther, 
A localization algorithm for $D$-modules.
Journal of Symbolic Computation {\bf 29} (2000), 721--728.
%
\bibitem{SST2}
M.~Saito, B.~Sturmfels, and N.~Takayama,
Hypergeometric polynomials and integer programing.
Compositio Mathematica {\bf 115} (1999), 185--204.
%
\bibitem{SST}
M.~Saito, B.~Sturmfels, and N.~Takayama,
{\it Gr\"obner Deformations of Hypergeometric Differential Equations},
Springer,
2000.
%
\bibitem{takayama-euler-darboux}
N.~Takayama,
Propagation of singularities of  solutions of the Euler-Darboux equation
and a global structure of the space of holonomic solutions I. 
Funkcialaj Ekvacioj {\bf 35} (1992), 343--403.
%
\bibitem{takayama-1992}
N.~Takayama,
An Approach to the Zero Recognition Problem by Buchberger Algorithms.
Journal of Symbolic Computation {\bf 14} (1992), 265--282.
%
\bibitem{WFS}
  The Wolfram Functions Site.\\
  {\tt http://functions.wolfram.com/08.05.26.0006.01}.
\end{thebibliography}
\end{document}